\documentclass[12pt]{amsart}

\usepackage{latexsym, amsthm, amsfonts, amsmath, amssymb}

\newtheorem{thm}{Theorem}
\newtheorem{prop}[thm]{Proposition}
\newtheorem{lem}[thm]{Lemma}

\theoremstyle{remark}
\newtheorem{example}[thm]{Example}

\theoremstyle{definition}
\newtheorem{defn}{Definition}

\DeclareMathOperator{\vol}{\mathrm{vol}}
\DeclareMathOperator{\diam}{\mathrm{diam}}

\title{On the asymptotic linking number}
\author{T.~Vogel}
\address{Mathematisches Institut, Universit\"at M\"un\-chen,
Theresienstr.~39, 80333~M\"un\-chen, Germany}
\email{tvogel@mpim-bonn.mpg.de}
\thanks{
        The author would like to thank Dieter Kotschick and Boris Khesin for their support.}

\begin{document}
\begin{abstract}
We prove a theorem formulated by V.~I.~Arnold concerning a relation between the asymptotic linking number and the Hopf invariant of divergence--free vector fields. Using a modified definition for the system of short paths, we prove their existence in the general case. 
\end{abstract}
\maketitle

\section{Introduction}
The aim of the present article is to complete a proof of a theorem stated by Arnold in~\cite{arn, ark}. Let $M$ be a closed oriented three--dimensional manifold having the real cohomology of a three--sphere, i.e. $H^1(M)=H^2(M)=0$. Let $\mu$ be a volume form on $M$. We consider two smooth vector fields $X, Y$ on $M$ which are divergence--free with respect to $\mu$, i.e. $di_X\mu=di_Y\mu=0$. Divergence--free vector fields arise as magnetic fields or as the velocity vector field of an incompressible fluid. This is actually the context where the Hopf invariant of divergence--free vector fields was considered for the first time, cf.~\cite{wol, mof}. It is defined as follows. Choose a 1-form $\eta_X$ with the property $d\eta_X=i_X\mu$. This is possible since $H^2(M)=0$. The Hopf invariant of $X, Y$ is given by the integral
\begin{equation*}
\int_M \eta_X\wedge i_Y\mu~.
\end{equation*} 

The asymptotic linking number can be described as follows. First, fix a set $\Sigma$ containing exactly one oriented path having starting point $p$ and end point $q$ for every pair $(p,q)$ of points of $M$. The set $\Sigma$ must have some additional properties, the proof that it is always possible to choose such a set $\Sigma$ is the main achievement of this paper. We thereby solve a problem which was apparently first noticed by P.~Laurence and that is discussed in~\cite{ark} on p.~146. The same problem arises in the work of Freedman and He, cf.~\cite{fh}, and is also solved by our construction.

For given periods $[0,T]$ respectively $[0,S]$ of time and starting points $x,y$ one considers the pieces of flow lines emerging from $x$ respectively $y$. After connecting the end point of each of these pieces with the starting point using the corresponding path in $\Sigma$, one almost always gets two disjoint closed oriented curves in $M$. In this case, they have a well defined linking number. The asymptotic linking number $\lambda_{X,Y}(x,y)$ is found by dividing the linking number of the two closed curves by $TS$ and letting $T$ and $S$ grow to infinity. The average linking number $\Lambda(X,Y)$ is defined to be the integral of $\lambda_{X,Y}(x,y)$ over $M\times M$.

The theorem formulated by Arnold yields an interpretation of the Hopf invariant of $X, Y$ in terms of the asymptotic linking number, more precisely, the Hopf invariant of $X, Y$ equals their average linking number. Recall that the classical Hopf invariant classifying continuous maps $S^3\rightarrow S^2$ can be interpreted as linking number, cf.~\cite{bt}. Suppose $f : S^3\rightarrow S^2$ is a differentiable representative of the homotopy class $[f]\in\pi_3(S^2)$ and $p,q\in S^2$ are regular values of $f$. Then the Hopf invariant of $f$ is equal to the linking number of $f^{-1}(p)$ and $f^{-1}(q)$.

The vector space of differential forms of degree $k$ will be denoted by $\Omega^k$, in case of compact support we will use $\Omega^k_c$~. Recall that if $N$ is an arbitrary manifold, a double form on $N\times M$ is a form on $M$ taking values in the forms on $N$. Here, one can interchange the roles of $N$ and $M$, a double form on $N\times M$ can also be considered as a form on $N$ taking values in the forms on $M$.

\section{Linking forms}

In this section, we define the linking number and the linking form. We give a new proof for the existence of a linking form. The proof in~\cite{ark} uses arguments of a more homological type and requires that $M$ has the integral homology of a 3-sphere. We use the theory of the Poisson equation on Riemannian manifolds.

The linking number of two closed oriented disjoint one-dimensional submanifolds $N_1, N_2\subset M$ can be defined as follows. For $i\in \{1,2\}$ choose disjoint tubular neighborhoods $W_i$ of $N_i$ in $M$ and representatives $\nu_i\in\Omega^2_c(W_i)$ of the Poincar\'e duals of $N_i\subset W_i$~. The extensions of $\nu_i$ to the whole of $M$ by zero will also be denoted by $\nu_i$~. Since $H^2(M)=0$, we can choose 1-forms $\gamma_i$ such that $d\gamma_i = \nu_i$~.
\begin{defn}          \label{d:lk}
The {\em linking number} of $N_1,N_2$ is denoted by $lk$ and defined by
\begin{equation*}
lk(N_1, N_2) = \int_M \gamma_1 \wedge \nu_2~.
\end{equation*}
\end{defn}
The linking number is independent of choices. The proof and the relation of Definition~\ref{d:lk} to other possible definitions of the linking number can be found in \cite{bt}.
\begin{defn}       \label{d:lform}
A {\em linking form} on $M$ is a double form $L$ on $M\times M$ such that for any two disjoint closed oriented one-dimensional submanifolds $N_1, N_2$ of $M$ the equality 
\begin{equation}              \label{g:lform1}
lk(N_1,N_2)=\int_{N_1} \int_{N_2} L
\end{equation}
holds.
\end{defn}
\begin{example}        \label{b:lfr3}
We want to describe a linking form on $\mathbb{R}^3$. For this example, let $\times$ be the usual cross-product of two vectors in $\mathbb{R}^3$, the standard inner product on $\mathbb{R}^3$ will be denoted by $\langle~,~\rangle$. If $V\in T_x\mathbb{R}^3$ and $W\in T_y\mathbb{R}^3$, we define $L_{\mathbb{R}^3}$ by 
\begin{equation}
L_{\mathbb{R}^3}(V,W)= \frac{1}{4\pi}\frac{\langle V,W\times(x-y)\rangle}{\|x-y\|^3}~.
\end{equation}
This double form is a linking form on $\mathbb{R}^3$, cf.~\cite{ark}.
\end{example}
We choose a metric on $M$ which we will use in several constructions. Let $*_y~,~\delta_y$ respectively $d_y$ denote the Hodge-$*$-operator, coderivative respectively derivative with respect to the second factor in $M\times M$. Let the double form $g$ be a fundamental solution for $\Delta : \Omega^*(M) \rightarrow \Omega^*(M)$. This means that for $\beta\in\Omega^*(M)$, the form $\zeta\in\Omega^*(M)$ defined by the integral
\begin{align*}
\zeta(x) = \int_{y\in M} \beta(y)\wedge*_yg(x,y)
\end{align*}
solves the Poisson equation $\Delta\zeta=\beta$. 
A fundamental solution exists if there are no nontrivial harmonic forms on $M$.  As the manifold $M$ has the real cohomology of a 3-sphere, it meets this condition by the Hodge theorem. For the explicit construction of $g$ cf.~\cite{der}. The solution of the equation $\Delta\zeta=\beta$ is uniquely determined if $\beta$ has degree $1$ or $2$. In particular, if $\beta\in\Omega^1(M)$ and $\zeta$ solves $\Delta\zeta=\beta$, then $d\zeta$ solves the Poisson equation with the inhomogeneous part $d\beta$. The equality 
\begin{equation}         \label{g:g1}
\int_{y\in M} d\beta(y)\wedge*_yg(x,y) = d\int_{y\in M} \beta(y)\wedge*_yg(x,y)
\end{equation}
follows. 
Both $g$ and $L=*_yd_yg$ are smooth away from the diagonal in $M\times M$. The singularities of $g$ respectively $L$ on the diagonal are of order $1$ respectively $2$, (cf. \cite{der}, p. 134). In particular both double forms are in $L^1$. 
The proof that $L$ is indeed a linking form relies on the following fundamental property of $L$.
\begin{lem}         \label{l:fundpropl}
For any 1-form $\omega$ on $M$ there is a function $h$ such that
\begin{equation}
\int_{y\in M} d\omega(y)\wedge *_yd_yg(x,y) = \omega(x) + dh(x)~.
\end{equation}
\end{lem}
\begin{proof}
Let $(M\times M)_{\varepsilon}$ denote the complement of an $\varepsilon$-neighborhood of the diagonal in $M\times M$. Let $\tilde{g}$ denote a smooth extension of $g\big|_{(M\times M)_\varepsilon}$ to the whole of $M\times M$. Since $\delta_y$ and $d_y$ are formally adjoint we have 
\begin{align*}
\int_{y \in M} d\omega(y) &\wedge *_y d_y g(x,y)
   =\int_{y \in M} d\omega \wedge *_y d_y \tilde{g}(x,y) + o(\varepsilon)\\
   &=\int_{y \in M} \delta d \omega(y) \wedge *_y \tilde{g}(x,y) +o(\varepsilon)\\
   &=\int_{y \in M} \Delta\omega \wedge *_yg(x,y) - \int_{y\in M} d\delta\omega(y)\wedge *_y\tilde{g}(x,y) + o(\varepsilon)\\ 
   &=\omega(x)-d\left(\int_{y\in M} \delta\omega(y)\wedge *_yg(x,y)\right) + o(\varepsilon)~.
\end{align*}
The last equality is due to \eqref{g:g1} and to the fact that $g$ is a fundamental solution of $\Delta$. Let $h$ be the integral in the second summand. Since $\varepsilon$ can be chosen to be arbitrarily small we have proved the lemma.
\end{proof}
\begin{thm}           \label{t:Llkf}
The form $L=*_yd_yg$ is a linking form on $M$.
\end{thm}
\begin{proof}
Let $\varepsilon>0$. For $i\in\{1,2\}$, let $N_i\subset M$ be disjoint one-dimensional submanifolds. We denote the $\varepsilon-$ball in $\mathbb{R}^2$ by $D^2_\varepsilon$. Choose disjoint tubular neighborhoods $W_i\cong N_i\times D^2_\varepsilon$ of $N_i$ such that the diffeomorphism is given by the geodesic exponential map restricted to the normal bundle of $N_i$ in $M$. Coordinates on $W_1$ respectively $W_2$ compatible with this product decomposition will be denoted by $(p,a)$ respectively $(q,b)$. As described in \cite{bt}, we choose representatives of the Poincar\'e duals $\omega_i\in\Omega^2_c(W_i)$ of $N_i$ such that for every $p\in W_i$
\begin{equation}    \label{g:Lp1}
\int_{\{p\}\times D^2_\varepsilon}\exp^*\omega_i=1~.
\end{equation}
Since $H^2(M)=0$ we can choose a form $\eta_1\in\Omega^1(M)$ satisfying  $d\eta_1=\omega_1$~. Then by Lemma~\ref{l:fundpropl} there is a function $h$ with the property
\begin{align*}
lk(N_1,N_2)=&\int_M \eta_1\wedge\omega_2 \\
           =&\int_{x\in W_1}\left(\int_{y\in W_2} \omega_1(y)\wedge*_yd_yg(x,y)-dh(x)\right) \wedge\omega_2(x)~.
\end{align*}
The summand containing $dh$ is zero by the Stokes's theorem. Since $W_1$ and $W_2$ are disjoint, we can treat $g$ like a smooth form in this integral. The last expression is therefore equal to
\begin{equation*}
\int_{(p,a)\in W_1}\int_{(q,b)\in W_2} \omega_1((p,a))\wedge\omega_2((q,b))\wedge(*_yd_yg((p,0),(q,0))+o(\varepsilon))~.
\end{equation*}
Because of the special choice~\eqref{g:Lp1} of the Poincar\'e duals 
\begin{equation*}
lk(N_1,N_2)=\int_{p\in N_1} \int_{q\in N_2} *_yd_yg + o(\varepsilon)
\end{equation*}
follows. This completes the proof because $\varepsilon>0$ can be arbitrarily small.
\end{proof}
\section{Systems of short paths}
For defining the asymptotic linking number, we need a specific way to get closed curves from pieces of flow lines of $X$ or $Y$. This will be achieved by closing up these pieces by prescribed paths. The paths used for this purpose must have certain properties described in Definition~\ref{d:sps}. 

As this definition is technical we first explain the meaning of some of the requirements. Integration over a path in $\Sigma$ is possible by (ii). The continuity condition in (iii) ensures that all the integrals considered are measurable functions on $M\times M$, the conditions \eqref{g:1}, \eqref{g:2} and \eqref{g:3} guarantee that the system of short paths does not contribute to the asymptotic linking number. Finally, property (iv) guarantees that the short paths do not intersect either the flow lines of $X$ and $Y$ or the other short paths too often. 
\begin{defn}        \label{d:sps}
A set $\Sigma$ of paths on $M$ is a {\em system of short paths} if it has the following properties.
\begin{itemize}
\item[(i)] For any two points $p,q\in M$ there is a unique path $\sigma(p,q)\in\Sigma$ starting at $p$ and ending at $q$. 
\item[(ii)] Each path in $\Sigma$ is piecewise differentiable.
\item[(iii)] The paths depend continuously on their endpoints almost everywhere. The limits 
\begin{align}
&\lim_{T,S \to \infty} \frac{1}{TS} \int_{\phi_{[0,T]}x} \int_{\sigma(\psi_Sy,y)} |L|= 0 \label{g:1} \\
&\lim_{T,S \to \infty} \frac{1}{TS} \int_{\sigma(\phi_Tx,x)} \int_{\psi_{[0,S]}y} |L|= 0 \label{g:2} \\
&\lim_{T,S \to \infty} \frac{1}{TS} \int_{\sigma(\phi_Tx,x)} \int_{\sigma(\psi_Sy,y)}  |L|= 0 \label{g:3}
\end{align}
exist in the $L^1$-sense.
\item[(iv)] The sets 
\begin{align*}
\mathcal{S}_{X,\Sigma} & = \left\{ (x,y) \in M\times M \big| \phi_{[0,T]}x \cap \sigma(\psi_Sy,y) \not= \emptyset \right\} \\
\mathcal{S}_{\Sigma,Y} & = \left\{ (x,y) \in M\times M \big| \sigma(\phi_Tx,x) \cap \psi_{[0,S]}y \not= \emptyset \right\} \\
\mathcal{S}_{\Sigma, \Sigma} & = \left\{ (x,y) \in M\times M \big| \sigma(\phi_Tx,x) \cap \sigma(\psi_Sy,y) \not= \emptyset \right\}
\end{align*}
have measure zero at any given time $T$ respectively $S$.
\end{itemize}
\end{defn}
Arnold's original definition requires that the integrals appearing in \eqref{g:1}, \eqref{g:2} and \eqref{g:3} are uniformly bounded for $T=S=1$ for almost all $x,y\in M\times M$. The intention of the following example is to illustrate a problem one encounters in this case.
\begin{example}
Let $\rho$ be a smooth function on $\mathbb{R}$ such that $\rho\equiv 0$ outside of $(1,3)$ and $\rho>0$ on $(1,3)$ such that $\rho(2)=1$. On $\mathbb{R}^3$ we use cylindrical coordinates $(r,\varphi,z)$. Consider the vector field
\begin{equation*}
X= \left\{ \begin{array}{ll} \rho(e^{2/z}r)e^{-1/z}\frac{\partial}{\partial\varphi} & \textrm{ if } z>0 \\
0 & \textrm{ if } z\le 0~.
           \end{array} \right.
\end{equation*}
It is smooth and divergence--free with respect to the standard volume form $rdr\wedge d\varphi\wedge dz$ on $\mathbb{R}^3$. At the origin it has a zero of infinite order. All flow lines of $X$ are periodic or degenerate. The period of orbits emerging from points with $z>0, r=2e^{-2/z}$ is $4\pi e^{-1/z}$. Hence for a smooth curve which is transversal to the plane $z=z_0>0$, the value of the integral in \eqref{g:1} for $T=S=1$ becomes arbitrarily large if one chooses $x=\left(2e^{-2/z_0},\varphi,z_0\right)$ and $z_0$ small enough, provided that the curve is close enough to the $z$-axes.   
\end{example}
Using our definition of a system of short paths we can prove their existence while the remaining parts of the construction remain unchanged. This closes a gap in the original proof of Arnold's theorem found in~\cite{arn, ark}. Our definition of the system of short paths has the additional advantage that the system can be chosen independently of $X,Y$, cf.~\cite{ark}.
\begin{thm}          \label{t:spsex}
Let $\Sigma$ be the set consisting of a geodesic of minimal length having starting point $p$ and end point $q$ for any $p,q\in M$.Then $\Sigma$ is a system of short paths. 
\end{thm} 
\begin{proof} 
The set of paths $\Sigma$ fulfills the conditions (i),(ii) and the continuity condition in (iii) in Definition~\ref{d:sps} by construction. Let $F_{T,S}$ denote the function 
\begin{equation*}     
\frac{1}{TS}\int_{\phi_{[0,T]}x}\int_{\sigma(y,\psi_Sy)}|L|
\end{equation*}
on $M\times M$. To prove that $\Sigma$ has the properties (iii) in Definition~\ref{d:sps}, we show first that $\lim_{T,S\to \infty}F_{T,S}= 0$ holds in the $L^1$-sense. This corresponds to~\eqref{g:1}. For $y\in M$ define  $S_1(y)=\left\{A\in T_yM \big| \|A\|\le1 \right\}$. Let 
\begin{equation}  \label{g:deff}
f(x,y)=\max_{A\in S_1(y)}\big\{\vert L(X(x),A)\vert \big\}~. 
\end{equation}
Let $r(x,y)$ denote the distance of two points $x,y\in M$. By~\cite{der} (p. 122 and 134),  there is a smooth double form $\alpha$ on $M\times M$ such that for $V\in T_xM$ and $A\in T_yM$ we have
\begin{equation*}
L(V,A)=\frac{\alpha(V,A)}{r(x,y)^2}~.
\end{equation*}
By the definition~\eqref{g:deff} of $f$ we have
\begin{align*}
f(x,y)=\frac{1}{r(x,y)^2}\max_{A\in S_1(y)}\big\{\vert\alpha(X(x),A)\vert\big\}
\end{align*}
and since $\alpha$ is a continuous double form, $r(x,y)^2f(x,y)$ is a continuous function on $M\times M$. Hence $\int_{x\in M} f(x,y)\mu$ is bounded above by a constant $C$ which does not depend on $y$. Let $\vol(M)$ denote the volume of $M$ with respect to the volume form $\mu$, the diameter of $M$ is denoted by $\diam(M)$. The geodesics are parameterized by $s\in[0,1]$, hence their velocity vector $\dot{\sigma}(\psi_Sy,y)(s)$ is no longer than $\diam(M)$. By Fubini's theorem we find 
\begin{align*}
\|F_{T,S}&\| = \int_{(x,y)\in M\times M} \frac{1}{TS} \int_{\phi_{[0,T]}x} \int_{\sigma(\psi_Sy,y)}|L|~\mu\wedge\mu \\
            &= \int_{(x,y)\in M\times M} \frac{1}{TS} \int_0^T \int_0^1 \left| L\big(X(\phi_tx),\dot{\sigma}(\psi_Sy,y)(s)\big) \right|~ds~dt~\mu\wedge\mu\\
            &\le \frac{1}{TS} \diam(M) \int_{y\in M} \int_0^1 \int_0^T \int_{x\in M} f\big(\phi_tx,\sigma(\psi_Sy,y)(s)\big)\mu~dt~ds~\mu\\
            &\le \frac{1}{TS}\diam(M)C\vol(M)T=\frac{1}{S}\diam(M)\vol(M)C~.
\end{align*}
The last expression becomes arbitrarily small if $S\to\infty$. This proves \eqref{g:1}. With an analogous argument one can show \eqref{g:2}. To prove \eqref{g:3} observe that the corresponding integral in Example~\ref{b:lfr3} remains bounded when two straight lines approach each other. The same is true for geodesics on the Riemannian manifold $M$. The set $\Sigma$ meets condition (iv) again by construction.
\end{proof}

The following theorem from ergodic theory will be needed. We only state a version sufficient in our situation, much more general statements are true (cf.~\cite{kre, tem}).
\begin{thm}                 \label{t:erg}
Let $\phi ,\psi$ be volume preserving flows on $M$ and $f$ an $L^1$ function on $M\times M$. Furthermore, let $(S_n), (T_n)$ be sequences of positive real numbers with $\lim_nT_n = \lim_nS_n=\infty$.
\begin{itemize}
\item[(i)] The limit
\begin{equation*}
\lim_{n\to \infty}\frac{1}{T_nS_n}\int_0^{T_n}\int_0^{S_n}f(\phi_tx,\psi_sy)~ds~dt
\end{equation*}
exists in $L^1$--sense. 
\item[(ii)]
Let $\tilde{f}$ be the limit function, then 
\begin{equation*}
\int_{M\times M}\tilde{f}~\mu\wedge\mu=\int_{M\times M}f~\mu\wedge\mu~.
\end{equation*}
\end{itemize}
\end{thm}
We omit the proof which is standard using methods from~\cite{geo} and~\cite{dus}. The Birkhoff ergodic theorem used in~\cite{arn, ark} implies that the limits exist in the pointwise sense. (Additionally, in the theorem of Birkhoff the sequences are required to be monotonous, cf.~\cite{tem}.) As the limits in Definition~\ref{d:sps} are supposed to exist in the $L^1-$ sense, Theorem~\ref{t:erg} is suitable.

\section{Asymptotic linking number and Arnolds theorem}
We now define the asymptotic linking number. To simplify notation we abbreviate $\sigma(\phi_Tx,x)$ by $\sigma_T$ and $\sigma(\psi_Sy,y)$ by $\sigma_S$. 
\begin{defn}        \label{d:aln}
The {\em asymptotic linking number} $\lambda_{X,Y}$ is the limit function
\begin{equation*}
\lambda_{X,Y}(x,y)=\lim_{T,S\to\infty}\frac{1}{TS}lk(\phi_{[0,T]}x\cup\sigma_T,\psi_{[0,S]}y\cup\sigma_S)~.
\end{equation*}
The {\em average linking number} $\Lambda(X,Y)$ is defined by
\begin{equation*}
\Lambda(X,Y)=\int_{M\times M} \lambda_{X,Y}~\mu\wedge\mu~.
\end{equation*}
\end{defn}
\begin{prop}
The limit in Definition~\ref{d:aln} exists in the $L^1$-sense and does not depend on the system of short paths. The limit function is integrable.
\end{prop}
\begin{proof}
By property (iv) of the short path system, the linking numbers are almost always well defined. Using the linking form and property (iii) in Definition~\ref{d:sps} we find
\begin{align}   \label{g:wohl1}
\lambda_{X,Y}(x,y) &= \lim_{T,S\to\infty} \frac{1}{TS} \int_{\phi_{[0,T]}x\cup\sigma_T} \int_{\psi_{[0,S]}y\cup\sigma_S}L\\
\nonumber & = \lim_{T,S\to\infty}\frac{1}{TS}\int_0^T\int_0^S L\big(X(\phi_tx),Y(\psi_sy)\big)~ds~dt~.
\end{align}
The last limit exists by Theorem~\ref{t:erg}. This theorem also implies that the limit function is $L^1$, just like the double form $L$ itself. Hence $\Lambda$ is well defined.
\end{proof}
We now give the definition of the Hopf invariant for divergence--free vector fields on $M$. Since the forms $i_X\mu$ and $i_Y\mu$ are closed and $M$ has the real cohomology of the three-sphere, these forms are exact. We can therefore choose 1-forms $\eta_X$ and $\eta_Y$ with the property that $d\eta_X = i_X\mu$ and $d\eta_Y=i_Y\mu$.
\begin{defn}           \label{d:hopfinv}
The {\em Hopf invariant} of $X,Y$ is defined to be the integral
\begin{equation*}
\mathcal{H}(X,Y)=\int_M \eta_X \wedge i_Y\mu~.
\end{equation*}
\end{defn}
\begin{prop}
The Hopf invariant of $X,Y$ does not depend on the choice of $\eta_X$. Furthermore $\mathcal{H}(X,Y)=\mathcal{H}(Y,X)$ and if $Z$ is another divergence--free vector field on $M$ we have $\mathcal{H}(X,Y+Z)=\mathcal{H}(X,Y)+\mathcal{H}(X,Z)$. 
\end{prop}
\begin{proof}
Let $\eta_X'$ be another 1-form with the property $d\eta_X'=d\eta_X=i_X\mu$. Hence $\eta_X-\eta_X'$ is closed. Since $H^1(M)=0$ , we can choose a function $f$ on $M$ such that $df=\eta_X - \eta_X'$. By the Stokes's theorem we find \begin{align*}
\int_M (\eta_X - \eta_X') \wedge i_Y\mu = \int_M df \wedge i_Y\mu = \int_M d(fi_Y\mu) = 0~.
\end{align*}
Hence $\mathcal{H}(X,Y)$ is well defined. Choose a 1-form $\eta_Y$ such that $d\eta_Y=i_Y\mu~.$ Since 
\begin{align*}
d(\eta_X \wedge \eta_Y) = i_X\mu \wedge \eta_Y - \eta_X \wedge i_Y\mu~,
\end{align*} 
symmetry follows from the Stokes's theorem and the fact that 2-forms commute with differential forms of any degree. The last statement follows from $i_{(Y+Z)}\mu = i_Y\mu + i_Z\mu$.
\end{proof}
Now we are in a position to state and prove Arnolds theorem. The proof is analogous to the proof found in~\cite{arn, ark}, but the asymptotic linking number is no longer defined as a pointwise limit but as a limit existing in the $L^1-$sense.
\begin{thm}                     \label{t:arn}
The average linking number of $X,Y$ is equal to the Hopf invariant of $X,Y$, i.e. 
\begin{equation*}
\Lambda(X,Y)=\mathcal{H}(X,Y)~.
\end{equation*}
\end{thm}
\begin{proof}
Recall that $d\eta_X=i_X\mu$ and $L=*_yd_yg$. Using Lemma~\ref{l:fundpropl} we find
\begin{align*}
\mathcal{H}(X,Y)&=\int_M \eta_X\wedge i_Y\mu\\
        &=\int_{x\in M} \left(\left(\int_{y\in M} i_X\mu(y)\wedge *_yd_yg(x,y) \right) -dh(x)\right)\wedge i_Y\mu~.
\end{align*}
Because $i_Y\mu$ is closed, the summand containing $dh$ does not contribute to the integral by the Stokes's theorem, hence
\begin{align}       \label{g:a1}
\mathcal{H}(X,Y) =\int_{x\in M}\int_{y\in M} L(X,Y)~\mu(x)\wedge\mu(y)~.
\end{align}
 By the second part of Theorem~\ref{t:erg}, the right hand side of~\eqref{g:a1} is equal to 
\begin{equation*}
\int_{(x,y)\in M\times M}\left( \lim_{T,S\to \infty}\frac{1}{TS}\int_0^T\int_0^S L\big(X(\phi_tx),Y(\psi_sy)\big)~ds~dt\right) \mu(x)\wedge\mu(y)
\end{equation*}
and this is equal to $\Lambda(X,Y)$ by~\eqref{g:wohl1}.
\end{proof}

\end{document}